\documentclass[11pt, a4]{amsart}
\usepackage[all]{xy}
\usepackage[a4paper, left=28mm, right=28mm, top=28mm, bottom=34mm]{geometry}
\usepackage{amsmath} 					
\usepackage{amsthm} 						
\usepackage{amssymb} 					
\usepackage{amscd} 						

\usepackage{epic, eepic}
\usepackage{graphicx} 		
\usepackage{here}					

\usepackage[all]{xy}							
\usepackage{longtable}

\usepackage{mathrsfs} 					
\usepackage{latexsym}						
\usepackage{type1ec} 
\usepackage[OT2,T1]{fontenc}
\usepackage[russian,english]{babel}

\usepackage{enumerate}					
\usepackage{verbatim}						
\AtBeginDocument{\shorthandoff{"}}			

\newcommand{\Rational}{\mathbb{Q}}					
\newcommand{\Integer}{\mathbb{Z}}					

\def\det{\mathop{\mathrm{det}}\nolimits}			

\def\Kernel{\mathop{\mathrm{Ker}}\nolimits}		
\def\sheafhom{\mathop{\mathscr{H}\kern -2pt om}\nolimits}		

\def\sheafend{\mathop{\mathscr{E}\kern -2pt nd}\nolimits}			

\def\sheafext{\mathop{\mathscr{E}\kern -2pt xt}\nolimits}			

\def\Left{\mathop{\mathrm{L} \kern -2pt}\nolimits}				
\def\Right{\mathop{\mathrm{R} \kern -2pt}\nolimits}			

\newcommand{\Cohomology}[2]{H^{#1}\! \left( {#2} \right)}

\newcommand{\strshf}{\mathcal{O}}				
\newcommand{\projsp}{\mathbb{P}}				

\def\sheafspec{\mathop{\mathscr{S}\kern -2pt pec}\nolimits}
\def\sheafproj{\mathop{\mathscr{P}\kern -2pt roj}\nolimits}

\def\Pic{\mathop{\mathrm{Pic}}\nolimits}			
\def\Jac{\mathop{\mathrm{Jac}}\nolimits}		

\def\GenLin{\mathop{\mathrm{GL}}\nolimits}			
\def\Matrix{\mathop{\mathrm{Mat}}\nolimits}		

\makeatletter
\def\iddots{\mathinner{\mkern1mu\raise\p@
    \hbox{.}\mkern2mu\raise4\p@\hbox{.}\mkern2mu
    \raise7\p@\vbox{\kern7\p@\hbox{.}}\mkern1mu}}
\def\adots{\mathinner{\mkern2mu\raise\p@\hbox{.} 
 \mkern2mu\raise4\p@\hbox{.}\mkern1mu
 \raise7\p@\vbox{\kern7\p@\hbox{.}}\mkern1mu}}
\makeatother

\newtheorem{thm}{Theorem}[section]			
\newtheorem{cor}[thm]{Corollary}					

\newtheorem{Alg}[thm]{Algorithm}			

\newtheorem{rmk}[thm]{Remark}					
\newtheorem{exa}[thm]{Example}					
	
\address{Department of Mathematics, Faculty of Science, Kyoto University, Kyoto 606-8502, Japan}
\title[Computation of linear determinantal representations]{An algorithm to obtain linear determinantal representations of smooth plane cubics over finite fields}
\date{\today}
\author{Yasuhiro Ishitsuka}
\email{yasu-ishi@math.kyoto-u.ac.jp}

\begin{document}

\maketitle
\begin{abstract}
We give a brief report on our computations 
of linear determinantal representations of smooth plane cubics
over finite fields. After recalling a classical interpretation of 
linear determinantal representations as rational points
on the affine part of Jacobian varieties, 
we give an algorithm to obtain all linear determinantal representations 
up to equivalence.
We also report our recent study on computations of 
linear determinantal representations of twisted Fermat cubics
defined over the field of rational numbers.
This paper is a summary of the author's talk at the JSIAM JANT workshop 
on algorithmic number theory in March, 2016. Details will appear elsewhere.
\end{abstract}
\section{Introduction}
Let $k$ be a field, and let
\begin{align*}
	F(X, Y, Z) &= a_{000}X^3 + a_{001}X^2Y + a_{002}X^2Z + \\
	&a_{011}XY^2 + a_{012}XYZ + a_{022}XZ^2 +
	a_{111}Y^3 + a_{112}Y^2Z + a_{122}YZ^2 + a_{222}Z^3
\end{align*}
be a {ternary cubic form} with coefficients in $k$
defining a smooth plane cubic $C \subset \projsp^2$.
The cubic $C$ is said to \textit{admit
a linear determinantal representation over $k$} if
there are a nonzero constant $0 \neq \lambda \in k$ and
three square matrices $M_0, M_1, M_2 \in \Matrix_3(k)$
of size 3 satisfying
\[
	F(X, Y, Z) = \lambda \cdot \det (M),
\]
where we put $M := XM_0 + YM_1 + ZM_2$.
Two linear determinantal representations $M, M'$ of $C$ are said to be 
\textit{equivalent} if there are invertible matrices $A, B \in \GenLin_3(k)$
satisfying
\(
	M' = AMB.
\)

Studying linear determinantal representations of smooth plane cubics is 
a classical topic in linear algebra and algebraic geometry 
(\hspace{-0.35pt}\cite{Vin89}, \cite{Dol12}).
Recently, they appear in the study of 
the derived category of smooth plane cubics
(\hspace{-0.35pt}\cite{Gal14}) and the theory of space-time codes (\hspace{-0.35pt}\cite{DG08}).
They have been studied from
arithmetic viewpoints (\hspace{-0.35pt}\cite{FN14}, \cite{II16}, \cite{Ish15}).

In this note, we study linear determinantal representations over finite fields. 
We give an algorithm to obtain all linear determinantal representations
of smooth plane cubics up to equivalence.
This paper is a summary of the author's talk at the JSIAM JANT workshop 
on algorithmic number theory in March, 2016. Details will appear elsewhere.

\section{Linear determinantal representations and rational points}\label{RatLDR}
Let $k$ be a field, and $F(X, Y, Z) \in k[X, Y, Z]$ a ternary cubic form
with coefficients in $k$ defining a smooth plane cubic $C \subset \projsp^2$. 
We fix projective coordinates $X, Y, Z$ of $\projsp^2$.
The following theorem gives an interpretation of 
linear determinantal representations 
of $C$ in terms of non-effective line bundles on $C$.
It is well-known at least when $k$ is 
an algebraically closed field of characteristic zero.
For the proof valid for arbitrary fields, see \cite[Proposition 3.1]{Bea00}, 
\cite[Proposition 2.2]{Ish15}.

\begin{thm}\label{Corr}
	There is a natural bijection between the following two sets:
	\begin{itemize}
		\item the set of equivalence classes of 
		linear determinantal representations of $C$ over $k$, and
		\item the set of isomorphism classes of
		non-effective line bundles on $C$ of degree $0$.
	\end{itemize}
\end{thm}

The bijection is obtained as follows:
we take a non-effective line bundle $\mathcal{L}$ of degree 0 on $C$.
Let $\iota \colon C \hookrightarrow \projsp^2$ be the given embedding.
We denote the homogeneous coordinate ring of $\projsp^2$ by
\begin{align*}
	R &:= \Gamma_*(\projsp^2, \strshf_{\projsp^2})\\
	&= \bigoplus_{n \in \Integer} \Cohomology{0}{\projsp^2,
	\strshf_{\projsp^2}(n)} \cong k[X,Y,Z].
\end{align*}
The graded $R$-module $N = \Gamma_*(\projsp^2, \iota_*\mathcal{L})
\cong \Gamma_*(C, \mathcal{L})$ has 
a minimal free resolution of the form
\begin{equation}\label{Eq: LDR}
	\xymatrix@=1.2em{
		0 \ar[r] & R(-2) \otimes_k W_1 \ar[r]^{\widetilde{M}} & 
		R(-1) \otimes_k W_0 \ar[r] & N \ar[r] & 0,
	}
\end{equation}
where $W_0, W_1$ are 3-dimensional $k$-vector spaces 
\cite[Proposition 3.1]{Bea00}.
The homomorphism $\widetilde{M}$ can be 
expressed by a square matrix $M$ of size $3$
with coefficients in $k$-linear forms in three variables $X, Y, Z$.
We can check $M$ gives a linear determinantal representation of $C$,
and its equivalence class depends only on 
the isomorphism class of the line bundle $\mathcal{L}$.

When $k$ is an algebraically closed field,
the set $\Pic^0(C)$ of isomorphism classes of line bundles 
on $C$ of degree 0 is parametrized by
the group $\Jac(C)(k)$ of $k$-rational points on 
the Jacobian variety $\Jac(C)$ of $C$,
and the only \textit{effective} line bundle of degree $0$ corresponds to 
the origin $\strshf$ of $\Jac(C)(k)$.
In general, there can be a difference between $\Pic^0(C)$ and
$\Jac(C)(k)$ which is measured by
the relative Brauer group (\hspace{-0.35pt}\cite[Theorem 2.1]{CK12}, 
\cite[Example 6.9]{Ish15}).
When $C$ has a $k$-rational point $P_0$, 
the relative Brauer group vanishes, 
and two sets $\Pic^0(C)$ and $\Jac(C)(k)$ are 
identified. We have a bijection 
\[
	C(k) \overset{1:1}{\to} \Pic^0(C)=\Jac(C)(k) \quad ; \quad 
	P \mapsto \strshf_C(P-P_0).
\]
Hence we obtain the following corollary.

\begin{cor}\label{Corr2}
	Let $C$ be a smooth plane cubic over $k$ 
	with a $k$-rational point $P_0 \in C(k)$.
	There is a natural bijection between the following two sets:
	\begin{itemize}
		\item the set of equivalence classes of 
		linear determinantal representations of $C$ over $k$, and
		\item the set $C(k) \setminus \{P_0\}$ 
		of $k$-rational points on $C$ different from $P_0$.
	\end{itemize}
\end{cor}

\section{An algorithm to obtain linear 
determinantal representations}

Let us make the bijection in Theorem \ref{Corr} explicit.
In this section, we give an algorithm 
to obtain linear determinantal representations
of smooth plane cubics over an arbitrary field $k$.
In this algorithm, we do not assume that $C$ has a $k$-rational point.
\begin{Alg}\label{Alg:LDR}\mbox{} \\[-13pt]
\begin{description}
	\item[Input:] a ternary cubic form $F(X,Y,Z)$ with coefficients in $k$
	defining a smooth plane cubic $C \subset \projsp^2$ 
	with respect to a fixed projective coordinates $X, Y, Z$, and
	a $k$-rational non-effective line bundle $\mathcal{L}$ on $C$ 
	of degree $0$.
	\item[Output:] a linear determinantal representation of $C$ over $k$ 
	corresponding to $\mathcal{L}$. 
	\begin{description}
		\item[Step 1 (Global Section)] Compute a $k$-basis $\{v_0, v_1, v_2\}$
		of the 3-dimensional $k$-vector space
		$\Cohomology{0}{C, \mathcal{L}(1)}$.
		\item[Step 2 (First Syzygy)] Compute a $k$-basis $\{e_0,$ $e_1,e_2\}$
		of the kernel of the multiplication map
		\begin{align*}
				\Cohomology{0}{C, \mathcal{L}(1)} \otimes
				& \Cohomology{0}{C, \strshf_C(1)} 
				\to
				\Cohomology{0}{C, \mathcal{L}(2)}.
		\end{align*}
		\item[Step 3 (Output Matrix)] 
		Write the $k$-basis $\{e_0, e_1, e_2\}$ as
		\[
			e_i = \sum_{j=0}^2 v_j \otimes l_{i, j}(X, Y, Z),
		\]
		where $l_{i,j}(X, Y, Z) \in \Cohomology{0}{C, \strshf_C(1)}$ 
		are $k$-linear forms.
		Output the matrix \[M = (l_{i,j}(X, Y, Z))_{0 \le i,j \le 2}.\]
	\end{description}
\end{description}
\end{Alg}
 By the sequence \eqref{Eq: LDR}, we have
 \[
 	W_0 \cong \Cohomology{0}{C, \mathcal{L}(1)}
 \]
 and
 \begin{align*}
 	W_1 \cong \Kernel \big( 
	\Cohomology{0}{C, \mathcal{L}(1)} \otimes
	& \Cohomology{0}{C, \strshf_C(1)} 
	\\ &\to
	\Cohomology{0}{C, \mathcal{L}(2)}
\big).
 \end{align*}
 Using $k$-bases of $W_0$ and $W_1$, we obtain an explicit matrix representation $M$ of
 the map $\widetilde{M}$ in \eqref{Eq: LDR}.
This $M$ gives a linear determinantal representation 
corresponding to $\mathcal{L}$ in the bijection of Theorem \ref{Corr}.

\section{An explicit formula on linear determinantal representations
of smooth plane cubics with rational points}
We apply Algorithm \ref{Alg:LDR} 
to a smooth plane cubic $C$ with a $k$-rational point $P_0$.
By changing the projective coordinates,
we may assume that $P_0 = [1:0:0]$ and 
the tangent line of $C$ at $P_0$ is the line $(Z=0)$.
The following theorem gives an explicit formula of the bijection in Corollary 
\ref{Corr2}. For the proof, see \cite[Theorem 4.1]{Ish16}.

\begin{thm}\label{Th: LDR2}
	Let $C \subset \projsp^2$ be a smooth plane cubic 
	over an arbitrary field $k$
	with a $k$-rational point $P_0 = [1:0:0]$.
	Assume that the tangent line of $C$ at $P_0$ is the line $(Z=0)$.
	We have the following formula for 
	a linear determinantal representation $M_P$ of $C$ over $k$
	corresponding to a point $P = [s:t:u] \in C(k) \setminus \{P_0\}$ 
	via Corollary \ref{Corr2}.
	\begin{description}
		\item[Case 1] 
		If $u \neq 0$, the equivalence class of 
		linear determinantal representations of $C$
		corresponding to $P$ is given by
		\begin{align}\label{DRformula1}
			M_P = \begin{pmatrix}
				0 & Z & -Y \\
				uY - tZ & 0 & L_0(X,Y,Z) \\
				u X - s Z & L_1(X, Y, Z) & L_2(X, Y, Z)
			\end{pmatrix},
		\end{align}
		where we denote
		\begin{align*}
			L_0(X, Y, Z) &:= -u^2 X -(a_{011}t^2 + a_{012}tu + a_{022}u^2 + su) Z,\\
			L_1(X, Y, Z) &:= u^2 a_{011} X + u^2 a_{111} Y+u (a_{111}t + a_{112}u) Z,\\
			L_2(X, Y, Z) &:= u(a_{011}t + a_{012}u) X + (a_{111}t^2 + a_{112}tu + a_{122}u^2)Z.
		\end{align*}
		\item[Case 2] 
		If $u=0$, the equivalence class of linear determinantal representations 
		of $C$ corresponding to $P$ is given by
		\begin{align}\label{DRformula2}
			M_P = \begin{pmatrix}
				0 & Z & -Y \\
				Z & a_{011}Y & \widetilde{L}_0(X, Y, Z)  \\
				\widetilde{L}_1(X,Y,Z) & \widetilde{L}_2(X, Y, Z) 
				& \widetilde{L}_3(X, Y, Z)
			\end{pmatrix},
		\end{align}
		where we denote
		\begin{align*}
			\widetilde{L}_0(X,Y,Z) &:= X + a_{012}Y + a_{022} Z, \\
			\widetilde{L}_1(X,Y,Z) &:= a_{011}X + a_{111}Y, \\
			\widetilde{L}_2(X, Y, Z) &:= a_{111}X 
			+(a_{012}a_{111} - a_{011}a_{112})Y,\\
			\widetilde{L}_3(X, Y, Z) &:= (a_{022}a_{111} - a_{011}a_{122})Y
			- a_{011}a_{222}Z.
		\end{align*}
	\end{description} 
\end{thm}

\begin{rmk}\label{Rm: Other1}
	Let $k$ be a field of characteristic not equal to 2 nor 3, and 
	\begin{align}\label{Eq: Wform}
		E \colon ( Y^2Z - X^3 - aXZ^2 - b Z^3 = 0)
		\subset \projsp^2
	\end{align}
	an elliptic curve over $k$ with origin $P_0 = [0:1:0]$
	defined by a Weierstrass equation.
	Let $P=[\lambda: \mu: 1] \in E(k) \setminus \{P_0\}$ be 
	a $k$-rational point on $E$.
	Galinat gave in \cite[Lemma 2.9]{Gal14}
	a representative of linear determinantal representations of $E$ over $k$
	corresponding to the divisor $P-P_0$ as
	\[
		M'_P := \begin{pmatrix}
			X - \lambda Z & 0 & -Y-\mu Z \\
			\mu Z - Y & X + \lambda Z & (a + \lambda^2) Z\\
			0 & Z & -X
		\end{pmatrix}.
	\]
	When $C$ has a $k$-rational flex, Theorem \ref{Th: LDR2} is equivalent
	to Galinat's formula. However, Theorem \ref{Th: LDR2} is also applicable 
	when $C$ has no $k$-rational flex.
	
	When $k$ is algebraically closed of characteristic not equal to 2 nor 3, 
	Vinnikov \cite{Vin89} gave other representatives.
\end{rmk}

\section{Applications to linear determinantal representations over finite fields}
Let $p$ be a prime number, and $m \ge 1$ a positive integer.
Let $\mathbb{F}_q$ be a finite field with $q = p^m$ elements.
It is well-known that any smooth plane cubic $C$ over $\mathbb{F}_q$ has 
an $\mathbb{F}_q$-rational point (\hspace{-0.35pt}\cite{Lan55}), hence
we can freely use Corollary \ref{Corr2} and Theorem \ref{Th: LDR2}
(at least after some changes of coordinates).

In \cite{Ish16}, we determine projective equivalence classes of 
smooth plane cubics over $\mathbb{F}_q$
with 0, 1 or 2 equivalence classes of linear determinantal representations.
We denote by $\mathrm{Cub}_q(n)$ the set of 
projective equivalence classes of smooth plane cubics over $\mathbb{F}_q$
with $n$ equivalence classes of linear determinantal representations.
By Corollary \ref{Corr2}, the number of elements $\#\mathrm{Cub}_q(n)$ 
coincides with the number of 
projective equivalence classes of smooth plane cubics over $\mathbb{F}_q$
with $n+1$ $\mathbb{F}_q$-rational points.
The latter can be determined by Schoof's formula \cite{Sch87}.
The following table summarizes the results of
our computations of $\#\mathrm{Cub}_q(n)$ when $0 \le n \le 2$. 
\vspace{-0pt}
\begin{table}[h]
	\caption{}
	\begin{center}
		\begin{tabular}{|c|c|c|c|c|c|c|}
			\hline
			 & $\mathbb{F}_2$ & $\mathbb{F}_3$ & 
			$\mathbb{F}_4$ & $\mathbb{F}_5$ & $\mathbb{F}_7$ & 
			$\mathbb{F}_q \; (q \ge 8)$
			\\ \hline
			$\#\mathrm{Cub}_q(0)$ & 1 & 1 & 1 & 0 & 0 & 0 
			\\ \hline
			$\#\mathrm{Cub}_q(1)$ & 1 & 1 & 1 & 1 & 0 & 0
			\\ \hline
			$\#\mathrm{Cub}_q(2)$ & 2 & 2 & 4 & 2 & 2 & 0
			\\ \hline
		\end{tabular}
		\label{Tb: FLDR1}
	\end{center}
\end{table}
\vspace{-10pt}

The following ternary cubic forms are representatives of 
$\mathrm{Cub}_q(0)$. 
They do not admit linear determinantal representations.
Each of them has only one rational point $[1:0:0]$.
\begin{itemize}
	\item $X^2Z+XZ^2+Y^3+Y^2Z+Z^3$ over $\mathbb{F}_2$.
	\item $X^2Z+Y^3-YZ^2+Z^3$ over $\mathbb{F}_3$.
	\item $X^2Z+XZ^2+Y^3+\omega Z^3$ over $\mathbb{F}_4$,
	where $\omega \in \mathbb{F}_4$ satisfies $\omega^2+\omega+1=0$.
\end{itemize}

The following ternary cubic forms are representatives of 
$\mathrm{Cub}_q(1)$. Each of them admits a unique equivalence class of
linear determinantal representations.
Their rational points are $[1:0:0]$ and $[0:0:1]$.

\begin{itemize}
	\item $X^2Z+XYZ+Y^3+Y^2Z+YZ^2$ over $\mathbb{F}_2$.
	\item $X^2Z-Y^3+Y^2Z+YZ^2$ over $\mathbb{F}_3$.
	\item $X^2Z+\omega XYZ+Y^3+Y^2Z+\omega YZ^2$ 
	over $\mathbb{F}_4 = \mathbb{F}_2[\omega]$.
	\item $X^2Z+Y^3+2YZ^2$ over $\mathbb{F}_5$.
\end{itemize}

The following ternary cubic forms are representatives of 
$\mathrm{Cub}_q(2)$. Each of them admits 
two equivalence classes of linear determinantal representations.

\begin{itemize}
	\item $X^2Z+XY^2+YZ^2$ over $\mathbb{F}_2$.
	The rational points are $[1:0:0], [0:1:0], [0:0:1]$.
	\item $X^2Z+XZ^2+Y^3$ over $\mathbb{F}_2$.
	The rational points are $[1:0:0], [1:0:1], [0:0:1]$.
	\item $X^2Z+XY^2+YZ^2+2XYZ$ over $\mathbb{F}_3$.
	The rational points are $[1:0:0], [0:1:0], [0:0:1]$.
	\item $X^2Z-XZ^2-XYZ-Y^3$ over $\mathbb{F}_3$.
	The rational points are $[1:0:0], [1:0:1], [0:0:1]$.
	\item $X^2Z+XY^2+\omega YZ^2$ over 
	$\mathbb{F}_4 = \mathbb{F}_2[\omega]$.
	The rational points are $[1:0:0], [0:1:0], [0:0:1]$.
	\item $X^2Z+XY^2+(\omega+1)YZ^2$ over 
	$\mathbb{F}_4 = \mathbb{F}_2[\omega]$.
	The rational points are $[1:0:0], [0:1:0], [0:0:1]$.
	\item $X^2Z+XZ^2+\omega Y^3$ over 
	$\mathbb{F}_4 = \mathbb{F}_2[\omega]$.
	The rational points are $[1:0:0], [1:0:1], [0:0:1]$.
	\item $X^2Z+XZ^2+(\omega + 1)Y^3$ over 
	$\mathbb{F}_4 = \mathbb{F}_2[\omega]$.
	The rational points are $[1:0:0], [1:0:1], [0:0:1]$.
	\item $X^2Z+XY^2+YZ^2-2XYZ$ over $\mathbb{F}_5$.
	The rational points are $[1:0:0], [0:1:0], [0:0:1]$.
	\item $X^2Z-XZ^2-2XYZ-Y^3$ over $\mathbb{F}_5$.
	The rational points are $[1:0:0], [1:0:1], [0:0:1]$.
	\item $X^2Z+XY^2+3YZ^2$ over $\mathbb{F}_7$.
	The rational points are $[1:0:0], [0:1:0], [0:0:1]$.
	\item $X^2Z-XZ^2+3Y^3$ over $\mathbb{F}_7$.
	The rational points are $[1:0:0], [1:0:1], [0:0:1]$.
\end{itemize}

We give some examples of linear determinantal representations
of the cubics in the above list without $k$-rational flexes.

\begin{exa}
	Consider the smooth plane cubic over $\mathbb{F}_2$ defined by
	\[
	 	X^2Z + XY^2 + YZ^2 = 0.
	\]
	This cubic has three $\mathbb{F}_2$-rational points;
	\[
		P_0 = [1:0:0], P_1 = [0:1:0], P_2 = [0:0:1].
	\]
	By Theorem \ref{Th: LDR2}, linear determinantal representations
	corresponding to $P_1, P_2$ are
	\[
		\begin{pmatrix}
			0 & Z & Y \\
			Z & Y & X \\
			X & 0 & Y
		\end{pmatrix},\\
		\begin{pmatrix}
			0 & Z & Y \\
			Y & 0 & X \\
			X & X & Z
		\end{pmatrix}.
	\]
\end{exa}
\begin{exa}
	Consider the smooth plane cubic over $\mathbb{F}_5$ defined by
	\[
	 	X^2Z + XY^2 + YZ^2 -2 XYZ = 0.
	\]
	This cubic has three $\mathbb{F}_5$-rational points;
	\[
		P_0 = [1:0:0], P_1 = [0:1:0], P_2 = [0:0:1].
	\]
	By Theorem \ref{Th: LDR2}, linear determinantal representations
	corresponding to $P_1, P_2$ are
	\[
		\begin{pmatrix}
			0 & Z & -Y \\
			Z & Y & X -2Y \\
			X & 0 & -Y
		\end{pmatrix}, 
		\begin{pmatrix}
			0 & Z & -Y \\
			Y & 0 & -X \\
			X & X & -2X + Z
		\end{pmatrix}.
	\]
\end{exa}

\section{twisted Fermat cubics over the field of rational numbers}
In this final section, we report our recent study on computations of 
linear determinantal representations of twisted Fermat cubics
defined over the field $\Rational$ of rational numbers.

Over the field $\Rational$ of rational numbers, some problems arise.
The main problem is that a line bundle $\mathcal{L}_P$ on $C$ is usually given by
the corresponding $k$-rational point $P$ on $\Jac(C)$, not on $C$.
This causes some problems in Step 1; the calculation of the $\Rational$-vector space
\(
	\Cohomology{0}{C, \mathcal{L}_P(1)}.
\)
Using the generalized Clifford algebra and the norm equation, 
we overcome these problems for twisted Fermat cubics.
We give simple examples. Details will appear elsewhere.
\begin{exa}
	Consider the smooth plane cubic over $\Rational$ defined by
	\[
		X^3 + Y^3 + Z^3 = 0.
	\]
	Its Jacobian variety is an elliptic curve whose Weierstrass equation is given by
	\begin{equation}\label{Jac1}
		Y^2Z - 9YZ^2 -X^3+27Z^3= 0
	\end{equation}
	(cf.~\cite{ARVT05}).
	Its $\Rational$-rational points are
	\[
		\strshf, [3 : 0:1], [3: 9:1]. 
	\]
	Let us take $P=[3: 0:1]$.
	The equivalence class of linear determinantal representations 
	corresponding to $P$ is represented by
	\[
		\begin{pmatrix}
			-X+2Y+Z & -2X+Y & X+Y \\
			X - Y & X+Z & -Y \\
			X & 2X+3Y & -2Y+Z
		\end{pmatrix}.
	\]
	The other point $[3: 9:1]$ corresponds to
	the transpose of the above matrix.
	This cubic does not admit a
	symmetric determinantal representation over $\Rational$ (cf.\ \cite{II16b}).
\end{exa}

\begin{exa}
	Consider the smooth plane cubic over $\Rational$ defined by
	\[
		2X^3 + 2^2Y^3 + Z^3 = 0.
	\]
	Its Jacobian variety is an elliptic curve whose Weierstrass equation is given by
	\[
		Y^2Z - 9 \cdot 2^3 YZ^2 - X^3 + 27 \cdot 2^6 Z^3 = 0
	\]
	(cf.~\cite{ARVT05}).
	It is isomorphic to \eqref{Jac1}, and its $\Rational$-rational points are
	\[
		\strshf, [3 \cdot 2^2 : 0:1], [3 \cdot 2^2: 9 \cdot 2^3:1]. 
	\]
	However, in \cite{Ish15}, we prove that these points do not correspond to
	linear determinantal representations over $\Rational$
	due to non-vanishing obstruction in the relative Brauer group of
	this cubic. This cubic does not admit a linear determinantal representation over $\Rational$.
\end{exa}

\begin{exa}
	Consider the smooth plane cubic over $\Rational$ defined by
	\[
		17X^3 + 17^2Y^3 + Z^3 = 0.
	\]
	Its Jacobian variety is an elliptic curve whose Weierstrass equation is given by
	\[
		Y^2Z - 9 \cdot 17^3 YZ^2 - X^3 + 27 \cdot 17^6Z^3 = 0
	\]
	(cf.~\cite{ARVT05}).
	It is isomorphic to \eqref{Jac1}, and its rational points are
	\[
		\strshf, [3 \cdot 17^2: 0:1], [3 \cdot 17^2: 9 \cdot 17^3:1]. 
	\]
	Let us take $P=[3 \cdot 17^2: 0:1]$.
	The equivalence class of linear determinantal representations 
	corresponding to $P$ is represented by
	\[
		\begin{pmatrix}
			3X-2Y+Z & -34X+153Y & 17X-51Y \\[3pt]
			\frac{1}{17}X - \frac{1}{17}Y & -3X-4Y+Z & 4X+7Y \\[3pt]
			\frac{1}{17}X + \frac{4}{17}Y & -2X-Y & 6Y+Z
		\end{pmatrix}.
	\]
	The other point $[3 \cdot 17^2 : 9 \cdot 17^3:1]$ corresponds to
	the transpose of the above matrix.
\end{exa}
%

\section*{Acknowledgements}
The author would like to thank sincerely to Professor Tetsushi Ito 
for various and inspiring comments. The work of the author was supported by 
JSPS KAKENHI Grant Number 16K17572.


\begin{thebibliography}{99}
	\bibitem[AR-VT05]{ARVT05}
	M.~Artin, F.~Rodriguez-Villegas and J.~Tate,
	On the Jacobians of plane cubics.
	Adv. Math., \textbf{198} (2005), Issue 1, pp.~366--382.

	\bibitem[Bea00]{Bea00}
	A.~Beauville,
	Determinantal hypersurfaces.
	Dedicated to W.\ Fulton on the occasion of 
	his 60th birthday.
	Michigan Math.\ J.\ \textbf{48} (2000), Issue 1, pp.~39--64.

	\bibitem[CK12]{CK12}
	M.~Ciperiani and D.~Krashen,
	Relative Brauer groups of genus 1 curves.
	 Israel J.\ Math., \textbf{192} (2012), pp.~921--949.
	  
	\bibitem[DG08]{DG08}
	A.~Deajim and D.~Grant,
	Space time codes and non-associative division algebras arising from elliptic curves. 
	Contemp.\ Math., \textbf{463} (2008), pp.~29--44.
	
	\bibitem[Dol12]{Dol12}
	I.~Dolgachev,
	Classical Algebraic Geometry: A Modern View.
	Cambridge University Press, Cambridge, 2012.
	
	\bibitem[FN14]{FN14}
	T.~Fisher and R.~Newton,
	Computing the Cassels--Tate pairing on the 3-Selmer
	group of an elliptic curve.
  Int.\ J.\ Number Theory. \textbf{10} (2014), No.~7, 
	pp.~1881--1907.
	
	\bibitem[Gal14]{Gal14}
	L.~Galinat,
	Orlov's Equivalence and Maximal Cohen--Macaulay
	Modules over the Cone of an Elliptic Curve.
	Math.\ Nachr.\ \textbf{287} (2014), No.~13, pp.~1438--1455.
	
	\bibitem[Ish15]{Ish15}
	Y.~Ishitsuka,
	A positive proportion of cubic curves over $\mathbb{Q}$ 
	admit linear determinantal representations.
	\textsf{arXiv:\ 1512.05167}, submitted.
	  
	\bibitem[Ish16]{Ish16}
	Y.~Ishitsuka,
	Linear determinantal representations of smooth plane cubics 
	over finite fields.
	\textsf{arXiv:\ 1604.00115}, submitted.  
	
	\bibitem[II16a]{II16}
	Y.~Ishitsuka and T.~Ito,
	The local-global principle for symmetric determinantal
	representations of smooth plane curves.
	Published online in The Ramanujan J.; 
	DOI: 10.1007/s11139-016-9775-3, 2016.
	
	\bibitem[II16b]{II16b}
	Y.~Ishitsuka and T.~Ito,
	On the symmetric determinantal representations of 
	the Fermat curves of prime degree.
	Int. J. of Number Theory, \textbf{12} (2016), Issue 4, pp.\ 955--967.
	
	\bibitem[Lan55]{Lan55}
	S.~Lang,
	Abelian varieties over finite fields.
	Proc.\ Natl.\ Acad.\ Sci.\ U.S.A.\ \textbf{41} (1955), No.~3, 
	pp.~174--176.
	
	\bibitem[Sch87]{Sch87}
	R.~Schoof,
	Nonsingular plane cubic curves over finite fields.
  J.\ of Comb.\ Theory, \textbf{46} (1987), Issue 2, pp.~183--211.
	
	\bibitem[Vin89]{Vin89}
	V.~Vinnikov,
	Complete description of determinantal representations of
	smooth irreducible curves.
	Linear Algebra Appl.\ \textbf{125} (1989), pp.~103--140.
	\end{thebibliography}
\end{document}